\documentclass[12pt]{article}
\usepackage{amsmath}
\usepackage{amssymb}
\usepackage{latexsym}
\usepackage{a4wide}
\input xy
\xyoption{all}

\providecommand{\bysame}{\leavevmode\hbox to3em{\hrulefill}\thinspace}
\def\lra{\longrightarrow}
\def\ra{\rightarrow}
\def\Omegabar{\overline{\Omega}}
\def\frak{\mathfrak}
\def\wt{\widetilde}

\def\SL{\mathrm{SL}}
\def\text{\textrm}
\def\ord{\mathrm{ord}}
\def\GL{\mathrm{GL}}

\def\Qbar{\overline{\Q}}

\def\Gal{\mathrm{Gal}}
\def\F{\mathbf F}
\def\Fbar{\overline{\F}}

\def\mod{\text{ mod }}

\def\rhobar{\overline{\rho}}
\def\qed{\hfill \square \ }

\def\Of{{\cal O}}
\def\C{\mathbb C}
\def\I{I}
\def\disc{\Delta}

\def\Q{\mathbf Q}
\def\Z{\mathbf Z}
\def\Zbar{\overline{\Z}}
\def\T{\mathbf T}
\def\twT{\widetilde{\T}}

\def\GL{\mathrm{GL}}

\def\mod{\text{ mod }}

\def\rat{\kappa}
\def\Lf{{\cal L}}
\def\E{\cal E}

\def\Div{\mathrm{Div}}
\def\Hom{\mathrm{Hom}}
\def\Pone{\mathbf{P}^1}
\def\new{\mathrm{new}}

\begin{document}
\newtheorem{sublemma}{Sub-lemma}
\newtheorem{theorem}{Theorem}[section]
\newtheorem{lemma}[theorem]{Lemma}
\newtheorem{df}[theorem]{Definition}
\newtheorem{cor}[theorem]{Corollary}
\newtheorem{conj}[theorem]{Conjecture}

\author{Frank Calegari \\ Matthew Emerton}
\title{The Hecke Algebra $\T_k$ has Large Index}
\maketitle
\abstract{  Let $\T_k(N)^{\new}$ denote the Hecke algebra acting
on newforms of weight $k$ and level $N$.
We prove that the power of $p$ dividing
the index of $\T_k(N)^{\new}$ 
inside its normalisation grows at least
linearly with $k$ (for fixed $N$), answering a question of Serre.
We also apply our method to give heuristic evidence towards
recent conjectures of Buzzard and Mazur.}

\section{Introduction.}

Let $\F$ be a finite field of characteristic $p$, and
let $\overline{\rho}: \Gal(\Qbar/\Q) \lra \GL_2(\F)$
be a modular Galois representation of tame level dividing $N$.
Let $X({\rhobar})$ denote the 
space of global Galois deformations unramified
outside primes dividing $Np$.  (More precisely, by $X(\rhobar)$ 
we mean the rigid analytic generic fibre of the corresponding
formal deformation space of $\overline{\rho}$.)
A fundamental question first raised explicitly by Mazur is to
understand the locus $\Omega$ of deformations inside 
$X({\rhobar})$ that are crystalline at $p$.
Assuming the Fontaine--Mazur conjecture \cite[Conj.~3c, p.~49]{fontaine-mazur},
$\Omega$ is precisely
the locus of modular points of level coprime to $p$.
The theory of the Eigencurve \cite{eigencurve}
provides a nice deformation
theory of such representations, but requires extra data:
in particular, the $U_p$ eigenvalue, or equivalently, a
Frobenius eigenvalue of the corresponding Dieudonn\'{e} module. The forgetful
map from the eigencurve to $X(\rhobar)$ is at most two to one
(see \cite[Thm.~6.11]{kisin}), but its
image is certainly complicated, consisting as it does
of an infinite union of modular arcs, each crossing to
form the ``infinite fern'' of Mazur and Gouv\^{e}a~\cite{GM}. 
A natural object of study
is the topological closure
$\Omegabar$ of $\Omega$ (in the space of all $\Qbar_p$-valued points
of $X(\rhobar)$, with its usual $p$-adic topology).
Two recent conjectures shed
light on the structure of $\Omegabar$.

For a deformation $\rho \in X(\rhobar)$, let $\Q_p(\rho)$ 
denote the field generated by the traces
$\mathrm{Tr}(\rho(\mathrm{Frob}_{\ell}))$,
for $\ell$ coprime to $Np$.

\begin{conj}[Buzzard \cite{Buzz}] \label{conj:buzzard} There exists a constant $c$
depending only on $N$ and $p$ such that for every
$\rho \in \Omega$ {\em (}or equivalently, every $\rho \in \Omegabar${\em )},
$$[\Q_p(\rho):\Q_p] < c.$$
\end{conj}

\begin{conj}[Mazur] \label{conj:mazur} Let $\rho \in \Omegabar$ be a classical
modular point. Then $\rho$ is tamely potentially semistable at $p$.
In other words, the local representation
attached to $\rho$ becomes semistable after restriction
to a tame extension of $\Q_p$.
\end{conj}

Conjecture~\ref{conj:buzzard} rules out the possibility that
$\Omegabar$ arises as the points of some deformation ring
$R$, since 
the rigid space associated
to $R$ would have to be at least one dimensional, and hence would contain
points defined over arbitrarily large extensions of $\Q_p$. 
This should be contrasted with the fact that
one does expect the set of crystalline representations
whose Hodge--Tate weights lie in some fixed
finite interval to be
parameterised by a
corresponding deformation space; indeed, assuming
the Fontaine--Mazur conjecture, the resulting set is finite
and the corresponding deformation space is given by a
product of twisted Hecke algebras.  Likewise, one should
contrast Conjecture~\ref{conj:mazur} with the results of
L.~Berger~\cite{Berger1} and
forthcoming work of Berger and P.~Colmez, in which it is shown
that the condition of being crystalline with bounded Hodge-Tate weights
is closed in local deformation space.
The subtlety of both conjectures thus lies
in the fact that while the level $N$ is fixed, the weights of
the modular forms that define the set $\Omega$ are arbitrary.

This paper proves some results which provide 
theoretical evidence for Conjectures~\ref{conj:buzzard}
and \ref{conj:mazur}, and should also be of
independent interest.

\emph{Acknowledgments:} The authors would like to thank the anonymous
referee for pointing out a stupidity in the proof of
Lemma~\ref{lemma:unproved}.

\section{Results}
Let $R$ be a finite reduced $\Z_p$-algebra.  We define
the \emph{index} of $R$ to be the index $[\tilde{R}:R],$
where $\tilde{R}$ is the normalisation of $R$.
We define the \emph{discriminant} of $R$ to be the
discriminant of the trace pairing on $R$ over $\Z_p$.
We will apply these notions when $R$ is one of various Hecke algebras.

For any
level $N\geq 1$ and
positive integer $k$, let $S_k(N)$ denote the space of
cusp forms of weight $k$
for the congruence subgroup $\Gamma_1(N)$ defined over $\Qbar_p$,
and let $S_k(N)^{\new}$ denote the subspace of $S_k$ spanned by the
newforms of conductor $N$.
Let $\T_k(N)$ denote the
finite $\Z_p$-algebra
$$\T_k(N):= \Z_p[T_2, \ldots ,T_n, \ldots] \subseteq \mathrm{End}_{\C}(S_k(N))$$
generated by Hecke operators,
and let $\T_k(N)^{\new}$ denote the quotient of $\T_k(N)$ that acts faithfully
on $S_k(N)^{\new}$.  The algebra
$\Q\otimes \T_k(N)^{\new} $ is reduced,
and has dimension $\dim_{\C} S_k(N)^{\new}$ over $\Q_p$.
We denote the discriminant of $\T_k(N)^{\new}$
by $\disc_k(N)^{\new}$, and its index by $\I_k(N)^{\new}$. 
What happens to the index as the weight approaches infinity?
\begin{theorem} \label{theorem:joch}
The $p$-adic valuation of the index grows
at least linearly in $k$. Equivalently,
$$\liminf_{k \rightarrow \infty} \frac{1}{k} \cdot  \ord_p(\I_k(N)^{\new})
\ > \ 0.$$
\end{theorem}

Let $\T_k(N)'$ denote the subring of $\T_k(N)$ generated by the Hecke operators
$T_n$ for $n$ prime to $N$.   This is a reduced subring
of $\T_k(N)$, and so we may similarly consider its discriminant
$\disc_k(N)'$ and its index $\I_k(N)'$.
For any divisor $M$ of $N$, the Hecke action on the
space $S_k(M)^{\new}$ (thought of as a subspace of $S_k(N)$) induces
a map $\T_k(N)' \rightarrow \T_k(M)^{\new}$,
and the product map $\T_k(N)' \rightarrow \prod_{M | N} \T_k(M)^{\new}$
is injective, and becomes an isomorphism after tensoring with $\Q$,  
and hence after passing to normalisations.
Thus $\I_k(N)' \geq \prod_{M|N} \I_k(M)^{\new}$,
and so Theorem~\ref{theorem:joch} implies
that $\ord_p(\I_k(N)')$ grows (at least) linearly in $k$.
This answers positively a question that
was first raised in print by Jochnowitz \cite{Joch}
(where she attributes
the question to Serre).
In that paper Jochnowitz also proves the weaker
statement that $\ord_p(\I_k(N)')$ becomes arbitrarily large as
$k$ approaches infinity.
Since $\disc_k(N)'$ (the discriminant
of $\T_k(N)'$) is always divisible by
$\I_k(N)'$, Theorem~\ref{theorem:joch} also has as a corollary that
$\disc_k(N)'$ grows (at least) linearly in $k$,
a result which was already known from the work
of Jochnowitz \cite{Joch}.



Our method of proving Theorem~\ref{theorem:joch}  is to first 
establish the following result:
\begin{theorem} \label{theorem:main}
Fix a level $N$, possibly divisible by $p$.
Let $f$ and $g$ be two normalised cuspidal $\Qbar_p$ eigenforms 
on $\Gamma_1(N)$ of arbitrary weight. Suppose moreover
that $f$ and $g$ are residually congruent; that is,
$f \equiv g \mod \frak p$, where $\frak p$ denotes the
maximal ideal of $\Zbar_p$.
Then there exists
a rational number $\rat > 0$ depending only on $N$ and $p$
such that
 for each $n$ prime to $p$,
$$a_{n}(f) \equiv a_{n}(g) \mod  p^{\rat} \Zbar_p.$$
\end{theorem}
The content of this theorem is the independence of
$\rat$ from the weights of $f$ and $g$. 
We view this result as some evidence towards Conjecture~\ref{conj:buzzard},
since it shows that 
at least modulo $p^{\kappa}$,
all eigenforms of fixed level have coefficients lying
in some bounded extension of~$\Q_p$. 

\

The constant $\rat$ in  
Theorem~\ref{theorem:main} can be made explicit, and
we use it to
produce the following evidence towards Conjecture~\ref{conj:mazur}:

\begin{theorem} \label{theorem:mazur} Assume the tame level $N = 1$.
Let $\rho:\Gal(\Qbar/\Q) \ra \GL_2(F)$ be
a modular point of $\Omegabar$ of weight $k$.
Assume moreover that $e(F) = 1$ \emph{(}i.e. $F/\Q_p$ is
unramified\emph{)}. Let $N(\rho)$ denote the conductor of $\rho$. If either
\begin{enumerate}
\item $p \ge 5$,
\item $p < 5$, $k = 2$, $N(\rho) \ne 27$,
\end{enumerate}
then $\rho$ is tamely potentially semistable.
\end{theorem}

The case $p \geq 5$ of this theorem is an immediate
consequence of Lemma~\ref{lemma:wild} below
(which when $F = \Q_p$ was proved as statement~(d) on p.~67 of
\cite{fontaine-mazur}).  It is the cases when $p =2$ or 3 that
are more delicate,
and which depend on certain
explicit versions of Theorem~\ref{theorem:main}.

\section{Eigenforms modulo powers of $p$}

The main difficultly in working with the space of
eigenforms over $\Zbar_p$ reduced 
modulo rational powers of $p$ is that this space
has no intrinsic geometric description.
The Serre--Deligne lemma guarantees that residual
eigenforms of weight at least $2$ automatically
lift to characteristic zero. This situation fails
miserably, however, for eigenforms over Artin rings.
For example, let $f$ and $g$ be two eigenforms
in $S_k(\Gamma_1(N),\Zbar_p)$
that are congruent modulo $p$ but distinct modulo $p^2$.
Then there exist  infinitely
many distinct eigenforms
$$(\alpha f + \beta g)/(\alpha + \beta) \in S_k(\Gamma_1(N),\Zbar_p/p^2),$$
only finitely many of which can lift to eigenforms of weight $k$ and
characteristic zero.
One way to avoid this problem is to not work directly with eigenforms
at all, and to instead work with the entire space $S_k$ of modular
forms. Any congruences satisfied by Hecke operators
on $S_k(\Gamma_1(N),\Zbar_p/p^{\kappa})$ will then automatically
be satisfied by eigenforms. This approach also has its
difficulties, however.

For example, when $N=1$ and $p=2$, it can be shown using
Coleman's theory \cite{coleman1} that
elements of $S_k(\SL_2(\Z),\Zbar_2/2)$ that are the reductions
mod $2$
of eigenforms in characteristic zero are always killed by $T_2$.
Let us see how close we can get to this congruence using
naive arguments:
It is easy to show that the Hecke operator
$T_2$ acts by zero on the space $S_k(\SL_2(\Z),\F_2).$
This implies in turn that $T_2$ induces a nilpotent operator
on $S_k(\SL_2(\Z),\Zbar_2/2).$  
However, the operator
$T_2$ on $S_k(\Zbar_2/2)$ is highly non-semisimple, and
the most that one can extract from this naive
analysis is that $T^n_2 = 0$, for some
$n$ that increases linearly with $k$. Correspondingly,
one can infer from this only that for eigenforms $f$,
$$a_2(f) \equiv  0 \mod 2^{1/n} \ \Zbar_2,$$
a congruence which is not independent of the weight of $f$,
and which is much weaker than the congruence cited at
the beginning of the paragraph.

\

In order to overcome the difficulties of the type that
occur in the naive argument of the preceding paragraph,
we adopt an approach that is
suggested by the arguments of Hatada \cite{Hatada}:
Namely, we choose a better $\T_k(N)$-invariant lattice 
$\Lambda \subset S_k(N)$,
for which the action of $\T_k(N)$ on $\Z_p/p \otimes \Lambda $ can be computed.
Ideally, the action of $\T_k(N)$ on $\Z_p/p \otimes \Lambda $ will
be semisimple, or as close to  this as possible.
It turns out that the space of
modular symbols is well adapted to our purpose, and it
is the method whereby Hatada~\cite{Hatada} obtains congruences
for modular forms of small level independent of the weight.
By formalizing aspects of this argument in terms of
geometry and
cohomology, we prove Theorem~\ref{theorem:main}.

\

\noindent \bf Proof of Theorem~\ref{theorem:main}\rm.
Without loss of generality we may assume that the congruence
subgroup $\Gamma_1(N)
\bigcap \Gamma(p)$ is torsion-free (by replacing $N$
by an appropriate multiple if necessary).
%
Let $Y(N,p)$ be the open modular curve that classifies
elliptic curves with a fixed point of order $N$ and
full level-$p$ structure.  We impose
no conditions on the value of the Weil pairing on the basis
elements giving the level-$p$ structure, and so 
$Y(N,p)$ is typically a disconnected curve, with $\phi(p)$
connected components. Concretely, $Y(N,p)$ is equal to the
quotient 
$$ \Gamma_1(N) \backslash \left( {\cal H} \times \GL_2(\Z/p) \right ),$$
where ${\cal H}$ denotes the upper half-plane, 
and $\Gamma_1(N)$ acts
on ${\cal H}$ through linear fractional transformations,
and on $\GL_2(\Z/p)$ by left multiplication.
Since $\Gamma_1(N) \bigcap \Gamma(p)$ is torsion-free,
$Y(N,p)$ does in fact represent the appropriate moduli problem,
and there exists a universal elliptic curve $\E$ over
$Y(N,p)$.
On $Y(N,p)$  we have the
standard rank two local system $\Lf$, corresponding to the
family of relative first cohomology groups of $\E$.
Let ${\Lf}_k$ denote the $k$th symmetric power
of $\Lf$; it is local system, free of rank $k+1$.
If $W_k$ denotes the $k$th symmetric power of the standard representation
of $\SL_2(\Z)$ on $\Z^2$, then ${\Lf}_k$ has the following concrete
description:
$${\Lf}_k := \Gamma_1(N) \backslash \left ( W_k \times {\cal H} \times
\GL_2(\Z/p) \right ).$$

Consider the short exact sequence of sheaves:
$$ 0 \longrightarrow {\Lf}_k \buildrel p \cdot \over \longrightarrow
{\Lf}_k \longrightarrow {\Lf}_k/p \longrightarrow 0.$$
Taking cohomology, and remembering that $H^2(Y(N,p),{\Lf}_k)$ vanishes 
(since its dual $H^0_c(Y(N,p),\check{{\Lf}}_k)$
vanishes (here $\check{{\Lf}}_k$ denotes the $\Z$-dual of
the local system of free $\Z$-modules ${\Lf}_k$),
$Y(N,p)$ being an open Riemann surface),
we obtain an isomorphism
$$H^1(Y(N,p), {\Lf}_k)/p \cong H^1(Y(N,p), {\Lf}_k/p).$$
This isomorphism is equivariant with respect to the prime-to-$p$
Hecke operators.

The local system ${\Lf}_k/p$ (the reduction of $\Lf$ modulo $p$)
is trivial over $Y(N,p)$, by definition of the
moduli problem that $Y(N,p)$ represents.
Thus 
$$H^1(Y(N,p), {\Lf}_k/p) \cong  W_k/p \otimes H^1(Y(N,p), \Z/p).$$
(The twisted coefficients ${\Lf}_k/p$ are actually untwisted, and so we
can pull them out of the cohomology.  The preceding concrete
description of ${\Lf}_k$ shows that the $k+1$-dimensional
$\Z/p$-vector space $W_k/p$ is the fibre of ${\Lf}_k/p$ over any point of 
$Y(N,p)$.)
Again,
this isomorphism is equivariant for all the prime-to-$p$ Hecke operators.
(Where on the right hand side, these ignore the $W_k/p$ factor, and
just act on the second factor.) In particular, if
$n$ is coprime to $p$, and $T_{n}$ is the $n$th
Hecke
operator acting on $H^1(Y(N,p), {\Lf}_k)/p$, then $T_n$
satisfies (modulo $p$) a polynomial of degree
independent of $k$, bounded explicitly by the dimension
of $H^1(Y(N,p), \Z/p)$. The Eichler--Shimura isomorphism
guarantees that any modular eigenform $f$
of weight $k$ and level $N$ corresponds to an 
eigenform in $\C \otimes H^1(Y(N,p),{\Lf}_k) $. In particular,
the eigenvalues of $T_{n}$ acting on  $f$ will
be among the eigenvalues of $T_{n}$ acting on
$H^1(Y(N,p),{\Lf}_k)$. Since $T_{n}$ (mod $p$) satisfies a
fixed polynomial independent of $k$, we infer that
there exists $\rat >0$ such that the reduction modulo $p^{\rat}$
of any eigenvalue of $T_{n}$ is determined explicitly
by its residue in $\Fbar_p$.
For example, one could take
$1/\rat$ to be the dimension of $H^1(Y(N,p),\Z/p)$.
This proves Theorem~\ref{theorem:main}. $\qed$

%

\

More precise values of $\rat$ can be extracted in particular
cases. For example:
\begin{lemma} \label{lemma:hatada} 
Assume that each connected component of $X(p)$ has genus $0$. Let
$f$ be a cuspidal eigenform 
of weight $k$ and level
$\Gamma(p)$ (not necessarily a newform). 
Then for all primes $\ell \equiv \pm 1 \mod p$, 
$$a_{\ell}(f) \equiv 1 + \ell \mod p \Zbar_p.$$
\end{lemma}

\begin{Proof}
If $X(p)$ has genus zero, then there is an isomorphism
$$H^1(Y(p),\Z/p) \simeq \wt{H}^0(C(p),\Z/p ),$$
where the target of this isomorphism denotes the reduced
$0$-dimensional cohomology of the set of cusps
$C(p) := X(p) \setminus Y(p)$. 
For $\ell \equiv \pm 1 \mod p$,
the action of $T_{\ell}$ on the cusps is given explicitly 
by $1 + \ell$. Thus $T_{\ell}$ satisfies the
polynomial $T_{\ell} - 1 - \ell = 0$ on
$$ W_k/p \otimes H^1(Y(p), \Z/p)$$
and the result follows. $\qed$ \end{Proof} 

\

Recall that the hypothesis of Lemma~\ref{lemma:hatada} is
satisfied precisely
for $p = 2,3,$ and $5$. 
Lemma~\ref{lemma:hatada} (phrased as a statement for these values of $p$)
was originally proved by Hatada~\cite{Hatada},
using a mixture of techniques. 
Our argument shows that Hatada's
explicit computations ultimately rely on the fact
that $X(p)$ has genus zero. However, our methods 
prove specific congruences for larger primes as well.
\begin{lemma} Let $p = 7$, and let 
$f$ be a cuspidal eigenform
of weight $k$ and level
$\Gamma(7)$ (not necessarily a newform).
Then for all primes $\ell \equiv \pm 1 \mod 7$,
$$a_{\ell}(f) \equiv 1 + \ell \mod \sqrt{7} \cdot \Zbar_7.$$
\end{lemma}
\begin{Proof} The curve $X(7)$ is the union of six
connected components each of genus $3$. There is
an exact sequence:
$$0 \ra H^1(X(7),\Z/7) \ra
H^1(Y(7),\Z/7) \ra \wt{H}^0(C(7),\Z/7 ) \ra 0.$$
The space $H^1(X(7),\Z/7)$ is ``accounted for'' by the
cubic twists of the weight $2$ 
form corresponding to the elliptic curve of 
conductor $49$ and CM by $\Z[(1+\sqrt{-7})/2]$.
Explicitly one finds that $T_{\ell}-1-\ell$ is
zero on $H^1(X(7),\Z/7)$ and on the cusps.
Thus $(T_{\ell}-1-\ell)^2 = 0$ on
$$ W_k/p \otimes H^1(Y(7), \Z/7)$$
and the result follows. $\qed$ \end{Proof}

\

 Serre conjectured \cite{Hatada}
that the congruence would be satisfied with $\sqrt{7}$
replaced by~$7$, although we are not able to prove this.

\

By passing to other curves of non-trivial genus,
and with a certain amount of non-trivial 
calculation
we are able to prove other
congruences, such as the following:

\begin{lemma} \label{lemma:unproved}
Let $f$ be a cuspidal eigenform
of weight $k$ and level $1$.
Then for all primes $\ell \equiv \pm 1 \mod 9$,
$$a_{\ell}(f) \equiv 1 + \ell \mod 9 \Zbar_3.$$
\end{lemma}

\begin{Proof}
Let us consider to begin with an arbitrary integer $N \geq 3;$
we will work on the modular curves $Y(N)$ and $Y(N^2)$
classifying elliptic curves with full level-$N$ structure
(respectively full level-$N^2$ structure).  We use notation
analogous to that introduced in the proof of Theorem~\ref{theorem:main}.

Consider the commutative diagram
\begin{equation}\label{eqn:diagram}
\xymatrix{H^1_c(Y(N^2), {\Lf}_k/N^2) \ar[d] \ar[r] \ar[dr] &
H^1(Y(N^2), {\Lf}_k/N^2)  \ar[d] \\
H^1_c(Y(N), {\Lf}_k/N^2) \ar[r] & H^1(Y(N), {\Lf}_k/N^2) ,}
\end{equation}
in which the vertical arrows exist because cohomology (with or
without compact supports) has
a covariant functoriality for proper maps (such as the map
$Y(N^2) \rightarrow Y(N)$),
and in which the diagonal arrow is defined as the composite of top and right-hand
(or equivalently, left-hand and bottom) arrows. 

\

{\em Claim.} For any local system ${\cal F}$ of finite $\Z$-modules
on $Y(N)$,
the push-forward
map $H^1(Y(N^2),{\cal F}) \rightarrow H^1(Y(N),{\cal F})$ is surjective.
(Here we have used ${\cal F}$ also to denote the pull-back of ${\cal F}$
to $Y(N^2)$.)

\smallskip

{\em Proof of claim.} A simple d\'evisage, using the fact
that $H^2(Y(N),{\cal F}) = 0$ (see the proof of Theorem~\ref{theorem:main}),
shows that it suffices to
prove the claim in the case when ${\cal F}$ is a local system
of $\F_{\ell}$-modules, for some prime $\ell$.
The map 
$H^1(Y(N^2),{\cal F}) \rightarrow H^1(Y(N),{\cal F})$ is then dual
to the pull-back map
$H^1_c(Y(N),\check{{\cal F}}) \rightarrow H^1_c(Y(N^2),\check{{\cal F}})$
(here ${\cal F}$ denotes $\F_{\ell}$-dual),
and so it suffices to show that this map is injective.
However, this latter map can be described in terms of modular symbols:
namely, if $\Div^0(\Pone(\Q))$ denotes the group of divisors of
degree zero supported on the elements of $\Pone(\Q)$, then
there is a commutative diagram
$$\xymatrix{
H^1_c(Y(N),\check{{\cal F}}) \ar[d]\ar[r]^-{\sim} &
\F_{\ell}[\pi_0(Y(N))]\,\check{} \otimes
\Hom_{\Gamma(N)}(\Div^0(\Pone(\Q)), \check{{\cal F}})\ar[d] \\
H^1(Y(N^2),\check{{\cal F}}) \ar[r]^-{\sim} &
\F_{\ell}[\pi_0(Y(N^2))]\,\check{} \otimes
\Hom_{\Gamma(N^2)}(\Div^0(\Pone(\Q)), \check{{\cal F}})}$$
in which the horizontal arrows are isomorphisms, and the right hand
vertical arrow is the tensor product of the dual of the surjection
$\F_{\ell}[\pi_0(Y(N^2))] \rightarrow
\F_{\ell}[\pi_0(Y(N))]$ (induced by the map $Y(N^2) \rightarrow
Y(N)$) and the inclusion
$\Hom_{\Gamma(N)}(\Div^0(\Pone(\Q)), \check{{\cal F}})
\subset
\Hom_{\Gamma(N^2)}(\Div^0(\Pone(\Q)), \check{{\cal F}}).$
In particular, it is injective, and thus
so is the left hand vertical arrow.

\

The preceding claim, applied to the local system ${\Lf}_k/N^2,$
shows that the right hand vertical arrow of diagram~\ref{eqn:diagram}
is surjective.

Let's make the following assumption:  

\

{\em Assumption}.  The diagonal arrow
$H^1_c(Y(N^2), {\Lf}_k/N^2) \rightarrow H^1(Y(N), {\Lf}_k/N^2)$ 
of diagram~\ref{eqn:diagram} vanishes.

\

{\em Claim.}  If the assumption holds, then for all primes
$\ell \equiv \pm 1 \pmod{N^2},$ the operator
$T_{\ell} - (1 + \ell)$ annihilates $H^1(Y(N), {\Lf}_k/N^2).$

\smallskip

{\em Proof of claim.}  The local system ${\Lf}_k/N^2$ is trivial on
$Y(N^2)$.  Thus we may rewrite the map
$H^1_c(Y(N^2), {\Lf}_k/N^2) \rightarrow H^1(Y(N^2), {\Lf}_k/N^2)$
as
$$W_k/N^2 \otimes H^1_c(Y(N^2), \Z/N^2) \rightarrow
W_k/N^2 \otimes H^1(Y(N^2), \Z/N^2).$$
Knowing that $T_{\ell} - (1 + \ell)$ kills
the cusps of $X(N^2)$, we see that its image on $H^1(Y(N^2), \Z/N^2)$
lies in the image of $H^1_c(Y(N^2), \Z/N^2)$ in $H^1(Y(N^2),\Z/N^2)$.
Thus the image of $T_{\ell} - (1 + \ell)$ on
$H^1(Y(N^2), {\Lf}_k/N^2)$ lies in the image of
$H^1_c(Y(N^2), {\Lf}_k/N^2)$ in $H^1(Y(N^2),{\Lf}_k/N^2)$. 
Since the right hand vertical arrow of diagram~\ref{eqn:diagram} 
is surjective, we find that the image of $T_{\ell} - (1 + \ell)$
on $H^1(Y(N), {\Lf}_k/N^2)$ is contained in the image
of $H^1_c(Y(N^2), {\Lf}_k/N^2)$ by the diagonal arrow of that diagram.
We are assuming that this vanishes, and so are done.

\

We turn to deriving sufficient conditions for the
assumption to hold.  
To this end, we suppose that each connected component of $X(N)$ has genus zero,
and consider the exact sequence of sheaves
$$0 \longrightarrow {\Lf}_k/N \buildrel N \cdot \over \longrightarrow
{\Lf}_k/N^2 \longrightarrow {\Lf}_k/N \longrightarrow 0.$$
Passing to cohomology, and also to cohomology with compact supports,
on $Y(N)$, we obtain the following diagram with exact rows:
$$\xymatrix{ 0 \ar[r] & H^1_c(Y(N),{\Lf}_k/N) \ar[r]\ar[d] &
H^1_c(Y(N),{\Lf}_k/N^2) \ar[r] \ar[d] & H^1_c(Y(N),{\Lf}_k/N) \ar[d] & \\
& H^1(Y(N),{\Lf}_k/N) \ar[r] &
H^1(Y(N),{\Lf}_k/N^2) \ar[r]  &
H^1(Y(N),{\Lf}_k/N) \ar[r] &  0.}$$ 
(As noted in the proof of Theorem~\ref{theorem:main},
since $Y(N)$ is an open Riemann surface,
$H^0_c$ and $H^2$ with coefficients in any local system always vanish.) 

Since ${\Lf}_k/N$ is trivial on $Y(N)$,
we have isomorphisms
$$H^1_c(Y(N),{\Lf}_k/N) \cong W_k/N \otimes H^1_c(Y(N),\Z/N)$$
and
$$H^1(Y(N),{\Lf}_k/N) \cong W_k/N \otimes H^1(Y(N),\Z/N).$$
Since the components of $X(N)$ have genus zero, the natural map
$$H^1_c(Y(N),\Z/N) \rightarrow H^1(Y(N),\Z/N)$$ vanishes,
and  so in the preceding diagram
the left-most and right-most vertical arrows both vanish.
Thus the middle arrow of that diagram factors through the natural map
$$H^1_c(Y(N), {\Lf}_k/N^2) \rightarrow  H^1_c(Y(N), {\Lf}_k/N),$$
and so the diagonal arrow of diagram~\ref{eqn:diagram} factors
through the composite
$$H^1_c(Y(N^2), {\Lf}_k/N^2) \rightarrow H^1_c(Y(N), {\Lf}_k/N^2) 
\rightarrow H^1_c(Y(N), {\Lf}_k/N).$$
Again, using the fact that ${\Lf}_k/N^2$ is trivial on $Y(N^2),$
and that ${\Lf}_k/N$ is trivial on $Y(N)$, we may rewrite this composite as
$$ W_k/N^2 \otimes H^1_c(Y(N^2), \Z/N^2) \rightarrow W_k/N \otimes
H^1_c(Y(N), \Z/N).$$
This map is obtained by tensoring through the natural map
$H^1_c(Y(N^2), \Z/N^2) \rightarrow H^1_c(Y(N), \Z/N)$ with
the projection $W_k/N^2 \rightarrow W_k/N$.
Thus to verify the assumption, it suffices to show that the former map
vanishes.

Exploiting the isomorphism (for a Riemann surface) between $H^1_c$
and $H_1$, we may rewrite this map as a direct sum over the connected
components of $X(N)$  of the maps
$$\Gamma(N^2)^{\rm ab}/N^2 \rightarrow \Gamma(N)^{\rm ab}/N.$$
Thus we see that our above assumption holds if the following two
conditions are met:  (i) Each component of $X(N)$ has genus zero; (ii) the image
of $\Gamma(N^2)^{\rm ab}$ in $\Gamma(N)^{\rm ab}$ lies in
$N \Gamma(N)^{\rm ab}.$

For any value of $N$, there is a natural isomorphism
$\Gamma(N)/\Gamma(N^2) \cong M_2(\Z/N)^0$ (the additive group of 
traceless $2 \times 2$ matrices).
Thus the commutator subgroup $\Gamma(N)^c$ of $\Gamma(N)$ is contained
in $\Gamma(N^2),$ and so we have the short exact sequence
$$0 \rightarrow \Gamma(N^2)/\Gamma(N)^c \rightarrow \Gamma(N)^{\rm ab}
\rightarrow M_2(\Z/N)^0 \rightarrow 0.$$

\

{\em Claim.}  If $N$ is such that 
$\Gamma(N)$ is free on three generators, then condition~(ii) above holds.

\smallskip

{\em Proof of claim.}
Since $\Gamma(N)$ is free on three generators, its abelianisation
is a free $\Z$-module of rank three, and so the kernel of the
surjection $\Gamma(N)^{\rm ab} \rightarrow M_2(\Z/N)^0$ must be
precisely $N \Gamma(N)^{\rm ab}$, since $M_2(\Z/N)^0 \simeq (\Z/N \Z)^3$. 

\

If we note that each component of $X(3)$ is of genus zero, and that
$\Gamma(3)$ is free on three generators, then we see that
the preceding claim proves the lemma.
$\qed$
\end{Proof}

\

This lemma allows one to prove (for example) that
the eigenforms of weight $2$ and level $243$
attached to elliptic curves 
cannot be approximated by $3$-adic eigenforms of level
$1$, since (in both cases) the coefficient $a_{19}$
fails to satisfy the required congruence. As
observed by Coleman and Stein~\cite{cstein}, the
congruences of Hatada can be used in a similar
way to eliminate the possible $2$-adic approximation
by forms of level one of the unique form of weight $2$ and level $32$.
In fact, by using  Lemma~\ref{lemma:unproved} and
Hatada~\cite{Hatada}, for $p=2$, $3$ 
one can show that 
all but one of the (finitely many) eigenforms of weight $2$, level
$p^n$ (for $n \ge 3$) with coefficients lying in an unramified
extension of $\Q_p$ cannot be $p$-adically
approximated by eigenforms of level one.
(The claim that there are finitely many such eigenforms
is justified in the remark following Lemma~4.2.) 
The exception is the unique form
of level $27$. This last example does not seem especially
anomalous; it could be dealt with by a slight
strengthening of Lemma~\ref{lemma:unproved}.
However,
just as our method fails to establish the conjecture of Serre mentioned above,
it also fails to prove the desired congruence. Thus,
although we are able to prove congruences for all primes $p$
and levels $N$, our method does not provide
a ``machine'' for proving  any particular congruence of this form.

\

\noindent \bf Proof of Theorem~\ref{theorem:joch}\rm.
Let us fix the level $N$ and the weight $k$, and also a maximal
ideal $\frak m$ of $\T_k(N)^{\new}$.  We will write
simply $\T$ to denote the localisation~$\T_k(N)^{\new}_{\frak m}$.

We write $\twT$ to denote the normalisation of $\T$; it is
a product of a finite number, say $d$, of finite DVR extensions of $\Z_p$.
Thus we may find an embedding of $\twT$ into $\Zbar_p^d$
such that the image of $\T$ is contained in
$\{(x_1,\ldots,x_d) \in \Zbar_p^d \, | \, x_1 \equiv \cdots \equiv x_d
\mod \frak p\}$,
where $\frak p$ denotes the maximal ideal of $\Zbar_p$.
Let $\T^{\circ}$ denote the subring of $\T$ generated by the prime-to-$p$
Hecke operators, let $\rat$ denote the rational number
of Theorem~\ref{theorem:main}, and let
$p^{\rat} \cap \twT$ denote the
intersection of $\twT$ with the ideal in $\Zbar_p^d$
generated by $p^{\rat}$.

\

{\em Claim:}  Let $\Of$ denote the $\Z_p$-subalgebra of $\Zbar_p$
obtained by projecting $\T^{\circ}$ onto the first factor of the product
$\Zbar_p^d$.
If we regard $\Of$ as being embedded diagonally into
the product $\Zbar_p^d$, then there is an inclusion
$\T^{\circ} \subseteq \Of + p^{\rat} \cap \twT.$

\smallskip

{\em Proof of claim:} Each of the projections of
$\T$ onto one of the factors of $\Zbar_p^d$ determines
a normalised Hecke eigenform, and by construction,
these eigenforms are all congruent modulo $\frak p$.
The claim thus follows from Theorem~\ref{theorem:main}.

\


The claim implies that
$$\T  = \T^{\circ}[T_p]  \subseteq \Of + p^{\rat} \cap \twT + \Of[T_p],$$
while elementary algebra implies that
$[\twT:\T]$ is divisible by
$[\twT/(p^{\kappa} \cap \twT):\T/(p^{\kappa} \cap \T)]$.
Yet
$$\T/(p^{\kappa}  \cap \T) \subseteq  \Of/(p^{\kappa} \cap \Of) +
\Of[T_p]/(p^{\kappa}\cap \Of[T_p]).$$
The first term is manifestly finite, and Theorem~\ref{theorem:main} shows
that it is independent of $k$.        
The theory of Coleman \cite{coleman1}, \cite{coleman2}
implies that the eigenforms of slope at most $\kappa$ fit
into finitely many analytic families. Thus the number of
such forms is bounded independently of the weight, and consequently
the second term is also finite and bounded
independently of $k$.

On the other
hand, if $\twT$ has rank $n$ over $\Z_p$, then a simple argument shows that
$\twT/(p^{\kappa} \cap \twT)$ has order at least $p^{\kappa n}$,
and thus
$\ord_p [\twT:\T] > \rat \cdot n - \epsilon$, 
for some explicitly computable constant $\epsilon$, independent of $k$.

Now take the product over all maximal ideals of $\T_k(N)^{\new}$.
From a well known result of Serre--Tate--Jochnowitz, 
the number of 
modular residual representations $\rhobar$ unramified outside $Np$ 
is finite, and this implies that the number of such maximal ideals
is bounded independently of the weight $k$.
On the other hand, the rank $n$ of $\T_k(N)^{\new}$ over $\Z_p$ grows
linearly with $k$,  and so the same is true of the
rank of its normalisation.  Thus the index
$I_k(N)^{\new}$ must grow at least linearly in $k$,
and Theorem~\ref{theorem:joch}\ is proved. $\qed$

\section{Crystalline Representations}

We  show in this section 
that there do not exist 
\emph{any} wildly potentially
semistable representations into $\GL_2(F)$ 
when $p \ge 5$, and $F$ is absolutely unramified.
Along with the discussion after Lemma~\ref{lemma:unproved}
this suffices to prove Theorem~\ref{theorem:mazur}.
This result and the  following lemma are undoubtedly known to the
experts (indeed, the case of the lemma when $F= \Q_p$ and $d = 2$
is proved as statement~(d) on p.~67 of \cite{fontaine-mazur}),
but we include a proof for lack of a 
reference in the literature.
We refer to \cite{fontaine1} 
for a discussion of semistable $p$-adic
representations of $\Gal(\Qbar_p/\Q_p)$ (and of the associated
terminology and notation), and to \cite{fontaine2} for a 
discussion of the Weil-Deligne group
representations attached to them.

\begin{lemma}\label{lemma:wild} Let $\rho: \Gal(\Qbar_p/\Q_p) \ra \GL_d(F)$
become semistable only after a wildly ramified extension.
Let $e$ denote the ramification index of $F$.
Then $p \le de + 1$.
\end{lemma}

\begin{Proof}  Denote the corresponding representation
by $V$. Let $K/\Q_p$ be a minimal Galois extension such
that the restriction of $V$ to $\Gal(\Qbar_p/K)$ is
semistable.
Let $K_0$ and $F_0$
denote the maximal unramified extensions of $\Q_p$ contained
in $K$ and $F$ respectively.
Without loss of generality, assume
that $K_0 \subseteq F_0$ (clearly this does not affect $e(F)$).
 Let $k_F$ be the residue
field of $K$, and $\Of = W(k_F)$  the ring of integers of
$F_0$.
 Let
$$D:=D_{st}(V|_K)
= (B_{st} \otimes_{\Q_p} V)^{G_K}.$$
$D$ is a free $K_0 \otimes_{\Q_p} F$ module of rank $d$.
By assumption $K_0 \subseteq F$, and thus there is
a natural isomorphism 
$$K_0 \otimes_{\Q_p} F \simeq \bigoplus_{\mathrm{Hom}(K_0,F)} F.$$
If elements of $\mathrm{Hom}(K_0,F)$
are denoted by $\sigma_i$, this
isomorphism is given explicitly by the
map $a \otimes b \mapsto 
[a \sigma_i(b)]_i$. If $e_{\sigma_j} = [\delta_{ij}]_i$ is
the $j$th idempotent, $D$ naturally decomposes as
a product
$$D \simeq \bigoplus_{\sigma_i} D_{\sigma_i}$$
where $D_{\sigma_i}:= e_{\sigma_i} D$.
Let $\sigma \in \mathrm{Hom}(K_0,F)$.
There is
a natural action of $\Gal(K/\Q_p)$ on $D$.
This action is $K_0$ semi-linear, and does not preserve
$D_{\sigma}$. Following \cite{fontaine2}, however, one may
adjust this action (using the crystalline Frobenius)
to obtain a \emph{linear} representation of the
Weil group $W_{\Q_p}$. In particular, the natural action
of $\Gal(K/K_0)$ (equivalently, the inertia
subgroup of $W_{\Q_p}$) on $D$ is linear and preserves 
$D_\sigma$. Since $D_{\sigma}$ is an $F$ vector space
of dimension $d$, we may also consider it as
an $F_0$ vector space of dimension $de$. 
Choosing a suitable lattice inside
$D_{\sigma}$ for the action of $\Gal(K/K_0)$ we
obtain a representation
$$\psi: \Gal(K/K_0) \rightarrow \GL_{de}(\Of).$$
If $\psi$ had non-trivial kernel, it
would correspond to some Galois subfield $L \subset K$
such that $\Gal(K/L)$ acted trivially on $D$. Yet 
then by descent $V$ would already be semistable over
$L$, contradicting the minimality assumption on $K$.
In particular, there must exist an element of
exact order $p$ inside
$$\GL_{d  e}(\Of).$$
The lemma then follows from the following well known
(and easy) result. $\qed$ \end{Proof}
\begin{lemma} \label{lemma:discrete} Let $\Of$ be a discrete valuation ring
with maximal ideal generated  by $p$. Suppose that $p > n+1$. Then
$\GL_n(\Of)$ has no elements of exact order $p$.
\end{lemma}

The preceding lemma is false if $p \leq n +1$.  Consequently,
if $p = 2$ or $3$, it is possible to have potentially semistable
representations $\Gal(\Qbar_p/\Q_p) \rightarrow \GL_2(\Q_p)$
that become semistable only after making a wildly ramified
extension. 
 However, a slightly more refined argument allows one to
bound the conductor of such a representation, and hence
to bound the power of $p$ that divides the level
of a newform whose Fourier coefficients lie in $\Q_p$.
The required generalisation of Lemma~\ref{lemma:discrete}
is that $\GL_n(\Of)$ contains (up to conjugacy)
a finite
number of finite $p$-groups. In particular,
$\rho$ becomes
semistable after some explicitly computable 
extension, from which one can explicitly
bound the conductor.
In particular, for any fixed weight,
there only finitely many
newforms of $p$-power level with coefficients in $\Q_p$.

\noindent \it Email addresses\rm:\tt \  fcale@math.harvard.edu

\hskip 22mm \tt \ emerton@math.northwestern.edu
\end{document}